\documentclass[12pt]{article}
\usepackage{amsmath,amsthm,amsfonts}

\numberwithin{equation}{section} \allowdisplaybreaks

\newtheorem{theorem}{\sc Theorem}[section]
\newtheorem{lemma}[theorem]{\sc Lemma}
\newtheorem{proposition}[theorem]{\sc Proposition}
\newtheorem{corollary}[theorem]{\sc Corollary}
\newtheorem{definition}[theorem]{\sc Definition}
\newtheorem{example}[theorem]{\sc Example}
\newtheorem{remark}[theorem]{\sc Remark}

\newcommand{\bet}{\begin{theorem}}
\newcommand{\eet}{\end{theorem}}
\newcommand{\blm}{\begin{lemma}}
\newcommand{\elm}{\end{lemma}}
\newcommand{\bprop}{\begin{proposition}}
\newcommand{\eprop}{\end{proposition}}
\newcommand{\bcor}{\begin{corollary}}
\newcommand{\ecor}{\end{corollary}}
\newcommand{\bdf}{\begin{definition}\rm}
\newcommand{\edf}{\end{definition}}
\newcommand{\bp}{\begin{proof}}
\newcommand{\ep}{\end{proof}}
\newcommand{\bex}{\begin{example}\rm}
\newcommand{\eex}{\end{example}}
\newcommand{\bremark}{\begin{remark}\rm}
\newcommand{\eremark}{\end{remark}}

\begin{document}

\title { B-Fredholm  Elements in Rings and Algebras}

\author{ M. Berkani}

\date{}

 \maketitle
\vspace{-12mm}

\begin{abstract}
In this paper, we  study  B-Fredholm elements in rings
and algebras. After  characterising these elements in terms of generalized Fredholm elements,
we will give sufficient conditions on a  unital primitive  Banach algebra $A$, under which we prove that   an element of $A$ is a B-Fredholm element  of index $0$ if and only if it is the sum of a Drazin
invertible element of $A$ and an element of the socle of $A$.

\end{abstract}

\renewcommand{\thefootnote}{}
\footnotetext{\hspace{-7pt}2010 {\em Mathematics Subject
Classification\/}:  Primary 47A53, 46H05.
\baselineskip=18pt\newline\indent {\em Key words and phrases\/}:
B-Fredholm, Banach algebra, index, inverse closed}

\section{Introduction}

This paper is a continuation of \cite{P42}, where we defined
B-Fredholm elements in  semi-prime Banach algebras, and we focused our
attention on the properties of the index. In particular,  we gave
a trace formula for the index of B-Fredholm operators. Here we
will consider in a first step  B-Fredholm elements in the case of
general rings, and then consider the case of primitive Banach algebras.

\noindent Let $X$ be a Banach space and let $L(X)$ be the Banach
algebra of Bounded linear operators acting on $X.$ In \cite{P7},
we have  introduced the class of linear bounded B-Fredholm
operators. If  $F_0(X)$ is the ideal of finite rank operators in $
L(X)$ and
 $\pi: L(X)\longrightarrow  A $ is the canonical projection,
where $ A= L(X)/F_0(X),$  it is well known by the Atkinson's
theorem \cite [Theorem 0.2.2,
p.4]{BMW}, that $T \in L(X) $ is a
Fredholm operator if and only if its projection
 $\pi(T)$ in the   algebra $ A $ is invertible.  Similarly, in the
 following result, we established an Atkinson-type theorem for B-Fredholm
operators.

\begin {theorem}\cite [Theorem 3.4]{P10}\label{thm1}: Let $T \in L(X).$ Then $T$ is a B-Fredholm operator
if and only if $\pi(T)$ is Drazin invertible in the algebra $
L(X)/F_0(X).$
\end{theorem}

\noindent Tacking into account this result and the definition of
Fredholm elements given in \cite{BA}, we defined in \cite{P42},
B-Fredholm elements in a semi-prime Banach algebra $A,$ modulo an ideal
$J$ of $A.$

\bdf \cite{P42}\label{def1} Let $A$ be  unital semi-prime Banach algebra and
let $J$ be an ideal of $A,$  and  $\pi: A\rightarrow A/J $ be  the
canonical projection. An element $a \in A$ is called  a B-Fredholm
element of  $A$ modulo the ideal $J$ if  its image $\pi(a)$ is
Drazin invertible in the quotient algebra $A/J.$ \edf

\noindent
Recall that  a ring $A$ is semi-prime if for $ u\in A, uxu= 0,$ for all $ x \in A$ implies that $u=0.$
A Banach algebra $A$ is called semi-prime if $A$ is also   a semi-prime ring.

\noindent  In a recent work \cite{CBS}, Cvetkovic and al., gave
in \cite[Definiton 2.3]{CBS} a definition of   B-Fredholm elements
in Banach algebras. However, their definition does not englobe the
class of B-Fredholm operators, since the algebra $L(X)/F_0(X)$ is
not a Banach algebra. That's why in our definition, we consider
general algebras, not necessarily being Banach algebras, so it
includes also the case of the algebra $L(X)/F_0(X).$

\noindent Recall that
Fredholm element in a semi-prime  ring $A$ were defined in \cite{BA}, as follows:
\bdf \cite[Definition 2.1]{BA}\label{def2} An element $a \in A$ is
said to be a Fredholm element of $A$ modulo $J$ if $ \pi(a)$ is
invertible in the  quotient ring $A/J,$ where
$\pi: A\longrightarrow  A/J $ is the canonical projection.\edf

\noindent Here we will use  Definition \ref{def2} and Definition \ref{def1} respectively, to define Fredholm elements and B-Fredholm elements in a unital ring .  As we will see in section 2, B-Fredholm elements in a
ring, and similarly to B-Fredholm operators acting on a Banach
space as observed in \cite[Proposition 3.3] {P12} , are related to
generalized Fredholm operators which had been defined  in \cite{CA} and
studied later in \cite{SC1}. Thus we will prove in section
2, that an element $a$ of a  ring  $A$ with a unit $e,$  is a B-Fredholm element of $A$ modulo
and ideal $J$ of $A$ if and only if there exists an integer $ n \in \mathbb{N},$ an
element $ c \in A$ such that  $a^nca^n-a^n \in J$ and $
e-a^nc-ca^n$ is a Fredholm element in $A$ modulo $J.$
Moreover, we will prove  a spectral mapping theorem for the Fredholm  spectrum and the B-Fredholm spectrum for  elements in  a unital Banach algebra.

\noindent In section 3, we will be concerned by  B-Fredholm
elements in a  unital primitive Banach algebra $A$  modulo the
socle of  $A.$ We  will give a condition on the socle of  $A$,
under which we prove that   an element of $A$ is a B-Fredholm
element  of index $0$ if and only if it is the sum of a Drazin
invertible element of $A$ and an element of the socle of $A,$
extending a similar decomposition given for B-Fredholm operators acting on a Banach space in \cite[Corollar 4.4]{P12}. Moreover, if   $p$ is any minimal idempotent in $A,$  for  $a
\in A,$ consider the operator $ \widehat{a}: Ap \rightarrow Ap,$
defined on the Banach space $Ap,$ by $\widehat{a}(y)= ay,$ for all
$ y \in Ap.$ Then,  we will give conditions under which there is an
equivalence between $a$ being a B-Bredholm element of $A,$ and
$\widehat{a},$ being a B-Fredholm operator on the Banach space
$Ap.$

\section {B-Fredholm elements in  ring }

Except when it is clearly specified,  in  all this section  $A$ will be  ring  with a unit $e,$
 $J$ an  ideal  of $A,$  and   $ \pi:
A\longrightarrow A/J$ will be the canonical projection.

 \bdf A non-empty subset ${\bf{R}}$ of   $A$
is called a regularity if it satisfies the following conditions:
\begin{itemize}
 \item If $ a \in A$ and  $n\geq 1$ is an integer,  then
  $ a \in {\bf R} $ if and only if $a^n \in {\bf R}, $

 \item  If $ a,b,c,d  \in \mathcal{A} $  are mutually commuting
  elements  satisfying  $ ac+bd = e, $ then $ ab\in {\bf R} $  if and
only if $a,b  \in {\bf R}. $

\end{itemize}

 \edf

\noindent Recall also that an element $ a \in A $ is said to be
Drazin invertible if there exists $b \in A $ and $ k \in
\mathbb{N}$ such that $ bab=b, ab=ba, a^kba=a^k.$

 \bet \label{thm2}The set of  Fredholm elements in $A$  modulo
$J$  is a regularity. \eet

\bp It is well known that the set of invertible elements in the
quotient ring $A/J$ is a regularity. Thus, its inverse image by
the ring homomorphism $ \pi$ is a regularity. \ep

\bet \label{thm3} The set of B-Fredholm elements in $A$ modulo $J$
 is  a regularity. \eet

\bp Similarly to  \cite[Theorem 2.3]{P10}, where it is proved that
the set of Drazin invertible elements in a unital  algebra is a
regularity, we can prove that in  the quotient ring $A/J$ the set
of Drazin invertible elements is a regularity. Thus, its inverse
image by the homomorphism $ \pi$ is also a regularity. \ep

\bprop  Let  $ a_1, a_2$ be B-Fredholm elements in $A$ modulo $J.$

i) If  $a_1a_2$ and $a_2a_1$ are elements of $J,$ then $ a_1 +
a_2$ is a B-Fredholm element in $A$ modulo $I$.

ii)  If $ a_1a_2 = a_2a_1,$ then $a_1a_2$ is a B-Fredholm element
in $A$ modulo $J.$

 iii) If $j$ is an element of $J,$  then $a_1+i$ is
a B-Fredholm element in $A$ modulo $J.$ \eprop

 \bp i) We have $\pi(a_1)\pi(a_2)= \pi(a_2)\pi(a_1)=0,$ From
\cite[Corollary 1]{DR},  it follows that $ \pi( a_1 + a_2)=
\pi(a_1) +\pi(a_2)$ is Drazin invertible in $A/J$. So  $ a_1 +
a_2$ is a B-Fredholm element in $A.$

ii)  We have $\pi(a_1a_2)= \pi(a_1)\pi(a_2)=\pi(a_2)\pi(a_1)$.
Similarly to \cite[Proposition 2.6]{P10},  it follows that
$\pi(a_1a_2)$ is Drazin invertible in $A/J.$ Hence $a_1a_2$ is a B-Fredholm element in $A$ modulo$J.$

 iii) If $i \in J,$  then
 $ \pi(a_1+i)= \pi(a_1)$. So $a_1+i$ is a B-Fredholm element in $A$ modulo$J.$

\ep

The following proposition is well known in the case Banach algebras. Its inclusion here witn a direct proof, in the case of unital rings, is  for a seek of completeness.

\bprop \label{prop1}

Let $A$ be a ring with a unit $e,$ and let $a \in A.$Then $a$ is
Drazin invertible in $A$ if and only if there exists an integer
$n \in \mathbb{N}^*,$ such that  $ A = a^nA \oplus N(a^n), $ where $ N(a^n)= \{ x
\in A \mid a^nx=0 \}.$ In this case there exists two idempotents $
p, q$ such that $ e= p+q , pq=qp=0$ and $A= pA \oplus qA.$

\eprop

\bp Assume that $a$ is Drazin invertible in $A.$ Then there exists
$b \in A $ and $ k \in \mathbb{N}$ such that $ bab=b, ab=ba,
a^kba=a^k.$ Without lose of generality, we can assume that $k=1.$
Let us show that $ A= aA \oplus N(a).$ Since $ aba=a$ then $ aA=
a^2 A.$ So if $x \in A,$ then $ ax= a^2t, $ with $ t\in A.$ Hence
$ a(x- at)=0,$ and $ x - at \in N(a).$ Therefore $ x= at +
(x-at).$ Moreover if $ x \in aA \cap N(a),$  then $ x=at,t \in A.$
Hence $ 0= bax= abat= at=x.$   Thus $ A = aA \oplus N(a).$ \ep

\noindent Conversely, assume that $ A= aA \oplus N(a).$ Then there
exists $ p \in aA,  q \in N(a),$ such that $ e=p+q.$ Then $p= p^2
+ qp,$ and $p-p^2= qp.$ As $ aA \cap N(a) = \{0\},$ then $p^2= p,$
and $qp=0.$ Similarly, we can show that $q^2=q,$ and $ pq=0.$
Moreover, if $ x\in A,$ then $x= ex= px+ qx.$ As $pA \cap qA= \{
0\},$ then $ A= pA \oplus qA.$ Thus, there exists $ r \in A,$ such
that  $a=pr,$ hence $pa= p^2r= pr=a.$ On the other side,we have $
ap= a(e-q)= a.$ Thus $ap=pa=a.$ Similarly, we have $ aq=qa=0.$
 Since $ A= aA \oplus N(a),$ then $aA= a^2A.$ So there
exists $ b \in aA,$ such that $ p= ab.$ Then $ a= pa= aba.$ We
have $ a(ba-p)= aba-ap=  0,$ so $ ba-p \in aA \cap N(a).$ Thus
$ba=p,$ and $ ba=ab=p.$ As $ bab- b= (ba-e)b= qb \in aA \cap
N(a),$ so $bab=b.$ Finally we have $aba=a, bab=b, ab=ba,$  and
$a$ is Drazin invertible in $A.$.

\bet \label{thm4} An element $a \in A$ is a B-Fredholm element of
$A$ modulo $J$ if and only if there exists an integer $ n \in
\mathbb{N}^*,$ an element $ c \in A$ such that  $a^nca^n-a^n \in J$
and $ e-a^nc-ca^n$ is a Fredholm element in $A$ modulo $J.$ \eet

\bp Assume that $ a$ is B-Fredholm element in $A$ modulo $J.$ Then
$\pi(a)$ is Drazin invertible in the quotient ring $A/J.$ Hence
there exists $b \in A $ and $ k \in \mathbb{N}^*$  such that $
\pi(b) \pi(a)\pi(b)= \pi(b) , \pi(a)\pi(b)= \pi(b) \pi(a),$ and
$\pi(a)^{k+1}\pi(b)= \pi(a)^k. $ So  $\pi(a)^k\pi(b)^k= \pi(b)^k
\pi(a)^k,  \pi(b)^k \pi(a)^k\pi(b)^k= \pi(b)^k, $ and
$\pi(a)^k\pi(b)^k \pi(a)^k= \pi(a)^k[\pi(b)\pi(a)]^k= \pi(a)^k. $
Let $c= b^k,$ then $\pi(e)- \pi(a)^k\pi(c)- \pi(c) \pi(a)^k=
\pi(e)$ is invertible in the quotient ring $A/J.$

\noindent Conversely suppose that there exists an integer $ n \in
\mathbb{N},$ and  an element $ c \in A$ such that
$\pi(a)^n\pi(c)\pi(a)^n=\pi(a)^n $ and $
\pi(e)-\pi(a)^n\pi(c)-\pi(c)\pi(a)^n$ is invertible in $A/J.$ Let
$ t= \pi(e)-\pi(a)^n\pi(c)-\pi(c)\pi(a)^n, s= t^{-1}$ and
let $ L_{\pi(a)^n}$ be the left multiplication in $A/J$ by $\pi(a)^n,
Im(L_{\pi(a)^n} )$ and $N(L_{\pi(a)^n})$ its image and kernel
respectively.
\noindent We have $ \pi(a)^n t= -\pi(a)^{2n} \pi(c),$ and $
\pi(a)^n= \pi(a)^n\pi(e)= \pi(a)^nt s= -\pi(a)^{2n} \pi(c)s.$
Hence $ \pi(a)^n A/J= \pi(a)^{2n} A/J,$  and so $Im(L_{\pi(a)^n}
)= Im(L_{\pi(a)^{2n}} ).$

 \noindent Similarly  we have
$t\pi(a)^n= -\pi(c)\pi(a)^{2n},$ and $ \pi(a)^n= \pi(e) \pi(a)^n=
st \pi(a)^n= -s\pi(c)\pi(a)^{2n}.$ Hence $N(L_{\pi(a)^n})=
N(L_{\pi(a)^{2n}}).$  Then it can be easily seen that
$ A/J= Im(L_{\pi(a)^n})
\oplus N(L_{\pi(a)^n}),$ where $\oplus $ stands for direct sum.
From Proposition \ref{prop1}, it follows that $ \pi(a)$ is Drazin
invertible in $A/J,$ and $ a$ is a B-Fredholm element in $A$
modulo $J.$ \ep

\vspace{2mm}

 \noindent Let us recall that an operator $T\in L(X)
$ has a generalized inverse if there is an operator $ S\in L(X) $
such that $TST=T.$
  In this case $S$ is  called a generalized inverse of $T$.
 It is well known that $T$ has a generalized inverse  if and only if
 $R(T)$ and  $N(T)$ are closed and complemented
subspaces of $X$. In \cite{CA}, S.R. Caradus  has defined the
following class of  operators :

\bdf  $ T\in L(X) $ is called a generalized Fredholm operator if
$T$ is relatively regular and there is a generalized
 inverse $ S $ of $T$ such that $ I-ST-TS$ is a Fredholm  operator.\\
\edf

 \vspace{-5 mm} \noindent In \cite{SC1} and \cite{SC2},
  this class of operators  had been  studied   and  it is  proved
\cite[Theorem1.1]{SC2} that an operator
   $ T\in L(X) $ is a generalized  Fredholm operator if and only if $T = Q \oplus F$,
  where  $Q$ is a finite dimensional nilpotent  operator and $F$ is  a Fredholm operator. Then Theorem \ref{thm4}, encourages us  to  consider
the following class of  elements in a  ring $A,$ with a unit $e.$

\bdf \label{}An element $a \in A$ is a generalized Fredholm element modulo
$J$ if there exists an element $b \in A$ such that $aba-a \in J$
and $ e-ab-ba$ is a Fredholm element in $A$ modulo $J.$ \edf

\noindent From Theorem \ref{thm4}, we obtain immediately the following
characterization of B-Fredholm elements

 \bet \label{thm5}  An element $ a
 \in A $ is a B-Fredholm element in $A$ modulo $J$ if and only if
there exists an integer $n \in \mathbb{N}^*$ such that $a^n$ is a generalized
Fredholm element in $A$ modulo $J.$
 \eet
\noindent Let $A$ be a  complex Banach algebra,  with unit $e,$ $a \in A,$ and let
$$\sigma_{F}(a)= \{ \lambda \in
\mathbb{C} \mid a- \lambda e \,\, \text {is not a Fredholm
element in } A \text \, {modulo} \, J\}, $$
\vspace{-5mm}
\text{and}

$$\sigma_{BF}(a)= \{ \lambda \in
\mathbb{C} \mid a- \lambda e \,\, \text {is not a B-Fredholm
element in } A \,\text {modulo} \, J\}, $$

 \noindent be respectively the Fredholm and the   B-Fredholm spectrum of $a.$  Then, we have the following result.

 \bet   Let $A$ be a unital  Banach algebra and $a \in A.$  If   $f$  an analytic function in
a neighborhood  of the usual spectrum $\sigma(a)$  of \,  $a$ which
 is non-constant on any connected component of the usual spectrum  $\sigma(a),$ of $a,$ then $
  f(\sigma_{{\bf BF}}(a))=\sigma_{\bf BF}(a)).$\\
 \eet
\vspace{ -5 mm} \noindent \bp From  Theorem  \ref{thm2} and Theorem \ref{thm3}, we know that the
set of Fredholm (resp; B-Fredholm) elements in $A$ is a regularity. Then the
corollary is a direct consequence of \cite[Theorem 1.4]{KM}. \ep

\section {B-Fredholm elements in primitive Banach algebras}

In this section, we will assume that $ A$ is a complex unital
primitive Banach algebra, with unit $e,$ and the ideal $J$ is
equal to its socle. Recall that an algebra is called primitive if
$ \{0 \}$ is a primitive ideal of $A.$ We will assume   that the
socle $J$ of $A$ is not reduced to $\{0\},$  so in this case and  from \cite{BMW},
$A$ possesses  minimal idempotents. A minimal idempotent $p$ of $A,$  is a non zero
idempotent $p$  such that $pAp= \mathbb{C}e.$  Recall also that it is
well know that a primitive Banach algebra is a semi-prime algebra.

\noindent Let  $p$ is any minimal idempotent in $A.$ For $a \in A,$ consider the operator $
\widehat{a}: Ap \rightarrow Ap,$ defined by $\widehat{a}(y)= ay,$ for all $ y \in Ap.$
We know from \cite[F.2.6]{BMW}, that if $ a$ is a Fredholm element in $A$, then $\widehat{a}$
is a Fredholm operator on the Banach space $Ap.$ However, the converse is in general false,
as shown in  \cite[F.4.2]{BMW}.

\noindent An  element $a \in A $ is said to be of finite rank if
the operator $ \widehat{a}$  is an operator of finite rank. We
know from \cite[Theorem F.2.4] {BMW}, that the socle of $A$ is
$ soc(A)= \{ x \in A\mid \widehat{x}$ {\text is of finite rank}\}.
The left regular  representation of the Banach  algebra $A$ on  the Banach
space $Ap$ is defined by  $\mathfrak{L}_r : A \rightarrow  L(Ap),$ such that
 $ \mathfrak{L}_r(x)=\widehat{ x}.$

\noindent For more details about the notions  from Fredholm
theory in Banach algebras used here , we refer the reader to \cite{BMW}.

\bdf \label{def3}\cite[2.1., p.283] {GR} Let $I$ be an ideal in a
Banach algebra $A.$ A function $ \tau : I\rightarrow \mathbb{C},$
is called a trace on $I$ if :

1)$ \tau(p)= 1$ if $p \in I,$ is an idempotent that is $p^2=p,$
and $p$ of rank one,

 2) $\tau (a+b)= \tau (a) + \tau (b), $ for all $ a, b \in I,$

3) $\tau(\alpha a)= \alpha \tau(a),$ forall $ \alpha \in
\mathbb{C}$ and $ a \in I,$

4) $\tau(ab)= \tau(ba),$ for all $ a \in I$ and $ b \in  A.$

\edf

\noindent From \cite [Section 3]{AM}, a trace function is defined
on the socle by: $ \tau(a)= \Sigma_{\lambda \in \sigma(a)}
m(\lambda, a)\lambda,$ for an element $a$ of the socle of $A,$
where $ \sigma(a)$ is the spectrum of $a,$ and $m(\lambda, a)$ is
the algebraic multiplicity of $\lambda$ for $a.$

\bdf \label{def31}  The index of a B-Fredholm element $a \in A$ is
defined by: $$ \mathbf{i}(a)= \tau( aa_0 - a_0a)= \tau ([a,
a_0]),$$

\noindent where $ a_0 $ is a Drazin inverse of $a$ modulo the
socle $J$ of $A.$\edf

\noindent From \cite[Theorem 2.3]{P42}, the index of a B-Fredholm element $a
\in A$ is well defined and is independant of the   Drazin inverse
$a_0 $  of $a$ modulo the ideal $J.$

\bdf Let $ a\in A$. Then $a$ is called a B-Weyl element
if it is a B-Fredholm element of index $0.$
\edf

\noindent The following theorem gives, under more hypothesis, a decomposition result for B-Fredholm elements of index $0$ in primitive Banach algebras, similar  to \cite[Theorem 3.1] {GR1} and  \cite[Theorem 3.2] {GR1}, given for Fredholm elements of index $0,$ in semi-simple Banach algebras. It  extends also a similar decomposition given for B-Fredholm operators acting on a Banach space in \cite[Corollar 4.4]{P12}.

\bet \label{thm3.1} Let $A$ be a unital primitive Banach algebra such that
 $ \mathfrak{L}_r(A)$ is Drazin inverse closed in  $L(Ap),$ and $\mathfrak{L}_r(soc(A))= F_0(Ap),$
where $F_0(Ap)$ is the ideal of finite rank operators in $L(Ap).$ Then an element
$a \in A,$ is a B-Weyl element if and only if $a=
b +c$ where $b$ is a Drazin invertible element of  $A$ and $c$ is
an  element of the socle  $J$ of $A.$ \eet

\bp Assume that $a$ is a  B-Fredholm element of index $0.$ Then from \cite[Lemma 3.2]{P42}, $
\widehat{a}$ is a B-Fredholm operator of index $0.$ From \cite[Corollary
4.4]{P12}, we have  $\hat{a}= S + F$, where S is a Drazin invertible
operator and F is a finite rank operator. Since  $A$ satisfies $\mathfrak{L}_r(soc(A))= F_0(Ap),$
 there exist $c \in J,$ such that $ F= \widehat{c}.$  As $S= \widehat{a-c}$ is Drazin invertible in $L(Ap)$, and $ \mathfrak{L}_r(A)$ is Drazin inverse closed in  $L(Ap),$ then $ \widehat{a-c}$ is
Drazin invertible in $ \mathfrak{L}_r(A)$. As the representation $\mathfrak{L}_r$ is faithful \cite[p. 30]{BMW}, then $ a-c$ is
Drazin invertible in $ A$.  Put $ b=a-c,$  then  $ a= b+ c,$ gives the
desired decomposition.

\noindent Conversely  if $a= b +c$ where $b$ is a Drazin
invertible element of  $A$ and $c$ is an  element of $J,$ then
from \cite[Proposition 3.3]{P42}, $a$ is a B-Fredholm element of
$A$ of index $0.$ \ep

\bex \label{exam1}
From \cite[Theorem F.4.3]{BMW}, if $A,$ is a  unital primitive C*-algebra, then
 $\mathfrak{L}_r(A)$ is inverse closed  in $L(Ap)$ and $\mathfrak{L}_r(soc(A))= F_0(Ap).$ Thus from \cite[Corollary 6]{RS},
$\mathfrak{L}_r(A)$ is Drazin inverse closed  in $L(Ap).$ Thus a primitive unital C*-algebra satisfies the hypothesis of Theorem \ref{thm3.1}
\eex

The aim of the  rest of this section , is to  establish a connection between B-Fredholmness of an element $a$ of $A$  and of the B-Fredholmness of the operator $ \mathfrak{L}_r(a)= \widehat{a}.$

\bet \label{thm3.2} Let $A$ be a primitive complex unital Banach algeba.
If $a$ is a B-Fredholm  element of $A$ modulo $J,$  then the
operator $ \widehat{a}$  is a B-Fredholm operator on the Banach
space $Ap.$ \eet

\bp If  $ a$ is a B-Fredholm element in $A$ modulo $J,$ then $a$
is Drazin invertible in $A$ modulo $J.$   From \cite[Theorem
F.2.4]{BMW}, we know that  J is exactly the set of elements $x$ of
$A$ such that $ \widehat{x}$ is an operator of finite rank. Then
$\widehat{a}$ is a Drazin invertible  operator modulo the ideal of
finite rank on $Ap.$ Thus from Theorem \ref{thm1}, $\widehat{a}:
Ap \rightarrow Ap $ is a B-Fredholm operator.
\ep
\vspace{3mm}

\noindent However, the converse of Theorem \ref{thm3.2} does not hold in general. To prove this, we use the same example as in \cite[Example F.4.2]{BMW}

\bex Let $T$ be the bilateral shift on the Hilbert space $l^2(\mathbb{Z}).$ Consider the closed unital subalgebra of $L(l^2(\mathbb{Z})),$ generated by $T$ and the ideal $K(l^2(\mathbb{Z}))$ of compact operators on $l^2(\mathbb{Z}).$ It follows from \cite[Example F.4.2]{BMW} that $A$ is a primitive Banach algebra and $ \{ \lambda \in \mathbb{C}: \mid \lambda \mid < 1\} \subset \sigma_A(T), $  the  spectrum of $T$ in $A.$ Hence $ 0 \in \sigma_A(T),$  and it is not an isolated point of $ \sigma_A(T).$ Therefore $a$ is not a B-Fredholm element of $A,$ otherwise and from \cite[Remark A, iii)]{P12}, if $ \lambda \neq 0$ and $\mid \lambda \mid $ is small enough, then $T- \lambda I$ is a Fredholm operator. But this impossible since from  \cite[Example F.4.2]{BMW},
$ \, \{ \lambda \in \mathbb{C}: \mid \lambda \mid  < 1\} \subset \sigma_{A/K(H)}(T+K(H)), $
 the spectrum of $T+K(H)$ in the Calkin algebra $ L(H)/K(H).$

\eex

In the following theorem, we give a necessary and sufficient condition, which ensures that  the converse of Theorem \ref{thm3.2} is true.

\bet \label{thm3.3} let $A$ be a primitive complex unital Banach algeba  satisfying $\mathfrak{L}_r(soc(A))= F_0(Ap).$
 Then the following two conditions are equivalent:

 i)For an element  $a \in A,$ if  $ \widehat{a}$  is a B-Fredholm operator on the Banach
space $Ap,$ then $a$ is a B-Fredholm element of $A.$

ii) Each element of the algebra
$\mathfrak{L}_r(A)/F_0(Ap)$ which is Drazin invertible in the algebra $L(Ap)/F_0(Ap),$
is also Drazin invertible in $\mathfrak{L}_r(A)/F_0(Ap).$
\eet

\bp It is clear that i) implies ii). So assume that for $ a \in A,$  $ \widehat{a}$  is a B-Fredholm operator on
the Banach space $Ap.$ From Theorem \ref{thm1},  $\widehat{a}$ is
Drazin invertible modulo the ideal of finite rank on $Ap.$     Since we assume that $ii)$ is true, then   $\widehat{a}$   is  Drazin invertible in $\mathfrak{L}_r(A)/ F_0(Ap).$  Thus   there exist $b \in A,$ such that $\widehat{a}\widehat{b}-\widehat{b}\widehat{a}, \widehat{b}\widehat{a}\widehat{b} - \widehat{b}, \widehat{a} ^{n+1}
\widehat{b}- \widehat{a}^n,$ are elements of $F_0(Ap).$
 As the representation $\pi$ is faithful, then
 $ab-ba, bab-b, a ^{n+1}b-a^n,$  are elements of $J.$  Thus $a$ is a B-Fredholm element of $A.$
 \ep
\bex  If $A$ is a  unital primitive $C^*$-algebra, then from \cite[Theorem F.4.3]{BMW}, we have $ \mathfrak{L}_r(soc(A)= F_0(Ap),$ where $p$ is a self-adjoint minimal idempotent of $A.$ Moreover if we assume that $A$ is commutative, and for  $a \in A,$  the operator $\widehat{a}$ is a B-Fredholm operator, then $\widehat{a}$ is a B-Fredholm multiplier on the $C^*-$ algebra $Ap.$  If $ \lambda \neq 0$ and $\mid \lambda \mid $ is small enough, then from \cite[Remark A, iii)]{P12},  $\widehat{a}- \lambda I$ is a Fredholm multiplier on the $C^*-$ algebra $Ap.$ From \cite[Corollary 5.105]{AI} , $\widehat{a}- \lambda I$ is of index is $0.$ Therefore and by \cite[Remark A, iii)]{P12}, $\widehat{a}$ is also of index $0.$ Hence by Example \ref{exam1} and Theorem \ref{thm3.1}, $a$ is a B-Weyl element of $A,$ and so a B-Fredholm element of $A.$
\eex

Similarly to Theorem \ref{thm3.3}, we have the following  result, which we give without proof.
\bet \label{thm3.4} let $A$ be a primitive complex unital Banach algeba  satisfying $\mathfrak{L}_r(soc(A))= F_0(Ap).$
 Then the following two  conditions are equivalent:

 i) For an element  $a \in A,$ if  $ \widehat{a}$  is a Fredholm operator on the Banach
space $Ap,$ then $a$ is a Fredholm element of $A.$

ii) Each element of the algebra
$\mathfrak{L}_r(A)/F_0(Ap)$ which is invertible in the algebra $L(Ap)/F_0(Ap),$
is also invertible in $\mathfrak{L}_r(A)/F_0(Ap).$
\eet

\bex  From \cite[Theorem F.4.3]{BMW}, if $A,$ is a  unital primitive C*-algebra, then
 $\mathfrak{L}_r(A)$ is inverse closed  in $L(Ap)$ and $\mathfrak{L}_r(soc(A))= F_0(Ap).$ Thus an element $a$ of $A$ is a Fredholm element if and only if $\widehat{a}$ is Fredholm operator.
\eex

 \baselineskip=12pt
\bigskip
\vspace{-5 mm }
 \baselineskip=12pt
\bigskip

{\small
\noindent Mohammed Berkani,\\
 \noindent Department of Mathematics,\\
 \noindent Science faculty of Oujda,\\
\noindent University Mohammed I,\\
\noindent Laboratory LAGA, \\
\noindent Morocco\\
\noindent berkanimo@aim.com,\\

\end{document}